\documentclass[11pt,twoside]{article}

\usepackage{pslatex}
\usepackage[T1]{fontenc}

\usepackage{amssymb}
\usepackage{amsmath}
\usepackage{bbm}
\usepackage{mathrsfs}
\usepackage{float}
\usepackage{xypic}

\makeatletter
\def\hsmash{\relax % \relax, in case this comes first in \halign
  \ifmmode\def\next{\mathpalette\mathhsm@sh}\else\let\next\makehsm@sh
  \fi\next}
\def\makehsm@sh#1{\setbox\z@\hbox{#1}\finhsm@sh}
\def\mathhsm@sh#1#2{\setbox\z@\hbox{$\m@th#1{#2}$}\finhsm@sh}
\def\finhsm@sh{\wd\z@\z@ \box\z@}
\makeatother

\makeatletter
\gdef\th@mychange{\normalfont\slshape
   \def\@begintheorem##1##2{\item
        [\hskip\labelsep \theorem@headerfont ##2. ##1  \,--\!--\!--\!--  ]}%
 \def\@opargbegintheorem##1##2##3{%
   \item[\hskip\labelsep \theorem@headerfont ##2. ##1\ {\upshape(}##3{\upshape)}. \,-----  ]}}
\makeatother

\RequirePackage{theorem}
\theoremstyle{mychange}

{\theorembodyfont{\rmfamily}\newtheorem{ttt}{}[section]}
{\theorembodyfont{\rmfamily}\newtheorem{nota}[ttt]{Notation.}}
{\theorembodyfont{\rmfamily}\newtheorem{defi}[ttt]{Definition.}}
{\theorembodyfont{\rmfamily}\newtheorem{remark}[ttt]{Remark.}}
{\theorembodyfont{\rmfamily}\newtheorem{rems}[ttt]{Remarks.}}
{\theorembodyfont{\rmfamily}\newtheorem{ex}[ttt]{Example.}}

{\theorembodyfont{\itshape}\newtheorem{fac}[ttt]{Fact.}}
{\theorembodyfont{\itshape}}
{\theorembodyfont{\itshape}\newtheorem{lem}[ttt]{Lemma.}}
{\theorembodyfont{\itshape}\newtheorem{prop}[ttt]{Proposition.}}
{\theorembodyfont{\itshape}\newtheorem{coro}[ttt]{Corollary.}}
{\theorembodyfont{\itshape}\newtheorem{theo}[ttt]{Theorem.}}

{\theorembodyfont{\rmfamily}}
{\theorembodyfont{\rmfamily}}
{\theorembodyfont{\rmfamily}\newtheorem{remo}[ttt]{Remark}}

{\theorembodyfont{\itshape}\newtheorem{obso}[ttt]{Observation}}
{\theorembodyfont{\itshape}\newtheorem{lemo}[ttt]{Lemma}}
{\theorembodyfont{\itshape}}
{\theorembodyfont{\itshape}\newtheorem{theoo}[ttt]{Theorem}}

\newcounter{abc}
\newenvironment{abc}{\begin{list}{\rm \alph{abc}) }{\usecounter{abc} \leftmargin=0.0pt \labelsep=0.0pt \listparindent=0.0pt \labelwidth=0.0pt \parsep=\smallskipamount \itemsep=0.0pt \topsep=0.0pt \partopsep=\smallskipamount}}{\end{list}}
\newcounter{iii}
\newenvironment{iii}{\begin{list}{\rm \roman{iii}) }{\usecounter{iii} \leftmargin=0.0pt \labelsep=0.0pt \listparindent=0.0pt \labelwidth=0.0pt \parsep=\smallskipamount \itemsep=0.0pt \topsep=0.0pt \partopsep=\smallskipamount}}{\end{list}}

\newcommand{\bP}{\mathop{\text{\bf P}}\nolimits}
\newcommand{\Spec}{\mathop{\text{\rm Spec}}\nolimits}
\newcommand{\rk}{\mathop{\text{\rm rk}}\nolimits}

\newcommand{\Gal}{\mathop{\text{\rm Gal}}\nolimits}
\newcommand{\Br}{\mathop{\text{\rm Br}}\nolimits}
\newcommand{\inv}{\mathop{\text{\rm inv}}\nolimits}
\newcommand{\Val}{\mathop{\text{\rm Val}}\nolimits}
\newcommand{\ev}{\mathop{\text{\rm ev}}\nolimits}
\newcommand{\Pic}{\mathop{\text{\rm Pic}}\nolimits}
\newcommand{\Div}{\mathop{\text{\rm Div}}\nolimits}
\renewcommand{\div}{\mathop{\text{\rm div}}\nolimits}

\newcommand{\Hom}{\mathop{\text{\rm Hom}}\nolimits}
\newcommand{\res}{\mathop{\text{\rm res}}\nolimits}
\newcommand{\cores}{\mathop{\text{\rm cores}}\nolimits}
\newcommand{\Frob}{\mathop{\text{\rm Frob}}\nolimits}
\newcommand{\cl}{\mathop{\text{\rm cl}}\nolimits}
\newcommand{\grad}{\mathop{\text{\rm grad}}\nolimits}
\newcommand{\vol}{\mathop{\text{\rm vol}}\nolimits}

\newcommand{\A}{\mathop{\text{\bf A}}}

\newcommand{\calD}{\mathscr{D}}
\newcommand{\calO}{\mathscr{O}}
\newcommand{\calQ}{\mathscr{Q}}
\newcommand{\calS}{\mathscr{S}}
\newcommand{\calU}{\mathscr{U}}

\newcommand{\bbA}{{\mathbbm A}}
\newcommand{\bbC}{{\mathbbm C}}
\newcommand{\bbF}{{\mathbbm F}}

\newcommand{\bbQ}{{\mathbbm Q}}
\newcommand{\bbR}{{\mathbbm R}}
\newcommand{\bbZ}{{\mathbbm Z}}

\renewcommand{\dagger}{\#}

\newcommand{\pmodulo}[1]{\nobreak\ifinner\mkern8mu\else\mkern18mu\fi
 (\textup{mod}\,\,#1)}

\newcommand{\br}{ }
\newcommand{\brr}{, }

\def\rightend#1#2{{%
 \leavevmode\nobreak\hskip .5em plus 1fil
 \penalty600 \hskip 0pt plus -1filll
 \vadjust{}\nobreak\hskip 0pt plus 1filll%
 #1\parfillskip=#2\relax \par}}

\def\eop{\ifmmode\rule[-22pt]{0pt}{1pt}\ifinner\tag*{$\square$}\else\eqno{\square}\fi\else\rightend{$\square$}{0pt}\fi}

\author{
\textbf{Andreas-Stephan Elsenhans}%
\footnote{E-mail: Stephan.Elsenhans@uni-bayreuth.de;
Website: {\scriptsize {\tt http://www.staff.uni-bayreuth.de/$\sim$btm216}}}~
and
\textbf{J\"org Jahnel}%
\footnote{E-mail: jahnel@mathematik.uni-siegen.de;
Website: {\scriptsize {\tt http://www.uni-math.gwdg.de/jahnel}}}}

\date{{\it${}^\dagger$Math.\ Inst., Universität Bayreuth, D-95440 Bayreuth, Germany\\
${}^\ddagger$FB6 Mathematik, Walter-Flex-Str.~3, D-57068 Siegen, Germany}\\
~\\
Received: July 7, 2009; Accepted: October 27, 2009}

\title{\textbf{\textsc{On the Brauer--Manin obstruction for cubic~surfaces}}\footnote{The computer part of this work was executed on the Sun Fire V20z Servers of the Gau\ss\ Laboratory for Scientific Computing at the G\"ottingen Mathematical Institute. Both authors are grateful to Prof.~Y.~Tschinkel for the permission to use these machines as well as to the system administrators for their~support.
}}

\begin{document}
\renewcommand{\thefootnote}{\fnsymbol{footnote}}

\maketitle

\begin{abstract}
We~describe a method to compute the Brauer-Manin~obstruction for smooth cubic surfaces
over~$\bbQ$
such that
$\Br(S)/\Br(\bbQ)$
is of order two or~four. This~covers the vast majority of the cases when this group is non-zero. Our~approach is to associate a Brauer class with every Galois invariant double-six. We~show that all order~two Brauer classes may be obtained in this~way. We~also recover Sir Peter Swinnerton-Dyer's result that
$\Br(S)/\Br(\bbQ)$
may take only five~values.% The~article is concluded by explicit~examples.
\end{abstract}


\section{Introduction}

\begin{ttt}
For~cubic surfaces, weak approximation and even the Hasse principle are not always~fulfilled. The~first example of a cubic surface violating the Hasse principle was constructed by Sir Peter Swinnerton-Dyer~\cite{SD1}. A~series of examples generalizing that of Swinnerton-Dyer is due to L.\,J.~Mordell~\cite{Mo}. An~example of a different sort was given by J.\,W.\,S.~Cassels and  M.\,J.\,T.~Guy~\cite{CG}.

A~way to explain these examples in a unified manner was provided by Yu.~I.~Manin in his book~\cite{Ma}. This~is what today is called the Brauer-Manin~obstruction. Manin's~idea is that a non-trivial Brauer class may be responsible for the failure of weak~approximation. We~will recall the Brauer-Manin~obstruction in some detail in section~\ref{zwei}.

An~important point is that only the factor
group~$\Br(S)/\Br(\bbQ)$
of the Grothendieck-Brauer~group of the cubic
surface~$S$
is~relevant. That~is isomorphic to the Galois cohomology
group~$H^1(\Gal(\overline{\bbQ}/\bbQ), \Pic(S_{\overline\bbQ}))$.
A~theorem of Sir Peter~Swinnerton-Dyer~\cite{SD} states that, for this group, there are only five~possibilities. It~may be isomorphic
to~$0$,
$\bbZ/2\bbZ$,
$\bbZ/2\bbZ \times \bbZ/2\bbZ$,
$\bbZ/3\bbZ$,
or~$\bbZ/3\bbZ \times \bbZ/3\bbZ$.
We~observed that, today, Swinnerton-Dyer's theorem from~1993 may easily be established by a script in~{\tt GAP}.

The~effect of the Brauer-Manin~obstruction has been studied by several~authors. For~example, for diagonal cubic surfaces, the computations were carried out by \mbox{J.-L.}~Colliot-Th\'el\`ene and his coworkers in~\cite{CKS}. In~this case,
$\Br(S)/\Br(\bbQ) = \bbZ/3\bbZ$.
The~same applies to the examples of Mordell or Cassels-Guy which were explained by the Brauer-Manin~obstruction in~\cite{Ma}.
\end{ttt}

\begin{ttt}
It~seems that, for~the cases that
$H^1(\Gal(\overline{\bbQ}/\bbQ), \Pic(S_{\overline\bbQ})) \cong \bbZ/2\bbZ$
or
$\bbZ/2\bbZ \times \bbZ/2\bbZ$,
no computations have been done up to~now. The goal of the present paper is to fill this~gap.

Our~starting point is a somewhat surprising~observation. It~turns out that
$\smash{H^1(\Gal(\overline{\bbQ}/\bbQ), \Pic(S_{\overline\bbQ}))}$
is of order two or four only in cases when,
on~$S$,
there is a Galois invariant \mbox{double-six}. This~reduces the possibilities for the action
of~$\Gal(\overline\bbQ/\bbQ)$
on the 27~lines. In~general, the automorphism group of the configuration of the 27~lines is the Weyl
group~$W(E_6)$~\cite[Theorem~23.9]{Ma}.
Among~the 350 conjugacy classes of subgroups
in~$W(E_6)$,
exactly 158 stabilize a double-six.

In~a previous paper~\cite{EJ}, we described a method to construct smooth cubic surfaces with a Galois invariant \mbox{double-six.} Our~method is based on the hexahedral form of L.~Cremona and Th.~Reye and an explicit Galois~descent. 
It~is able to produce examples for each of the 158 conjugacy~classes.

Among~them, however, there are 56 which even stabilize a~sixer. Those~may be constructed by blowing up six points
in~$\bP^2$
and, thus, certainly fulfill weak~approximation. There~are 26 further conjugacy classes which lead to
$H^1(\Gal(\overline{\bbQ}/\bbQ), \Pic(S_{\overline\bbQ})) = 0$.
\end{ttt}

\begin{ttt}
In~this article, we compute the Brauer-Manin obstruction for each of the 76~cases such that
$H^1(\Gal(\overline{\bbQ}/\bbQ), \Pic(S_{\overline\bbQ})) \cong \bbZ/2\bbZ$
or~$\bbZ/2\bbZ \times \bbZ/2\bbZ$.
We~start with two ``model cases'' for the Brauer groups
$\bbZ/2\bbZ$
and~$\bbZ/2\bbZ \times \bbZ/2\bbZ$.
These~are the maximal
subgroup~$U_1 \subset W(E_6)$
stabilizing a double-six and the maximal
subgroup~$U_3 \subset W(E_6)$
stabilizing a triple of azygetic double~sixes~\cite{Ko}.

In~both cases, we compute the Brauer group~explicitly. This~means, we produce representatives which we describe as Azumaya~algebras. We~then show that every
subgroup~$H \subset W(E_6)$
which leads to a Brauer group of order four is actually contained
in~$U_3$.
Recall~that every
subgroup~$H \subset W(E_6)$
leading to a Brauer group of order two is contained
in~$U_1$.
Finally,~we prove the main result that the restriction map is bijective in each of the~cases.
\end{ttt}

\begin{ttt}
The~article is concluded by examples showing the effect of the Brauer-Manin~obstruction. It~turns out that, unlike the situation described in~\cite{CKS} where a Brauer class of order three typically excludes two thirds of the adelic points, various fractions are~possible.
\end{ttt}

\section{The Brauer-Manin obstruction -- Generalities}
\label{zwei}

\begin{ttt}
For~cubic surfaces, all known counterexamples to the Hasse principle or weak approximation are explained by the following~observation.
\end{ttt}

\begin{defi}
Let~$X$
be a projective variety
over~$\bbQ$
and~$\Br(X)$
its Grothen\-dieck-Brauer~group. Then,~we will~call
$$\ev_\nu \colon \Br(X) \times X(\bbQ_\nu) \longrightarrow \bbQ/\bbZ \, , \quad
(\alpha, \xi) \mapsto \inv\nolimits_\nu (\alpha |_\xi)$$
the {\em local evaluation map\/}.
Here,~$\inv_\nu \colon \Br (\bbQ_\nu) \to \bbQ/\bbZ$
(and~$\inv_\infty \colon \Br (\bbR) \to \frac12\bbZ/\bbZ$)
denote the canonical~isomorphisms.
\end{defi}

\begin{obso}[{\rm Manin}{}]
Let\/~$\pi\colon X \to \Spec (\bbQ)$
be a projective variety
over\/~$\bbQ$.
Choose an element\/
$\alpha \in \Br (X)$.
Then,~every\/
$\bbQ$-rational
point\/
$x \in X(\bbQ)$
gives rise to an adelic point\/
$(x_\nu)_\nu \in X(\A_\bbQ)$
satisfying the~condition
$$\sum_{\nu\in\Val(\bbQ)} \!\!\!\! \ev_\nu (\alpha, {x_\nu}) = 0 \, .$$
\end{obso}

\begin{rems}
\label{standard}
\begin{iii}
\item
It~is obvious that altering
$\alpha \in \Br(X)$
by some Brauer~class
$\pi^*\rho$
for
$\rho \in \Br(\bbQ)$
does not change the obstruction defined
by~$\alpha$.
Consequently,~it is only the factor group
$\Br(X) / \pi^*\!\Br(\bbQ)$
which is relevant for the Brauer-Manin~obstruction.
\item
The~local evaluation map
$\ev_\nu \colon \Br(X) \times X(\bbQ_\nu) \to \bbQ/\bbZ$
is continuous in the second~variable.
\item
Further,~for every projective
variety~$X$
over~$\bbQ$
and
every~$\alpha \in \Br (X)$,
there~exists a finite
set~$S \subset \Val(\bbQ)$
such that
$\ev (\alpha, \xi) = 0$
for every
$\nu \not\in S$
and~$\xi \in X(\bbQ_\nu)$.

These~facts imply that the Brauer-Manin obstruction, if present, is an obstruction to the principle of weak~approximation.
\end{iii}
\end{rems}

\begin{lem}
\label{75}
Let~$\pi\colon S \to \Spec \bbQ$
be a non-singular cubic~surface. Then, there is a canonical~isomorphism
$$\delta\colon H^1(\Gal(\overline{\bbQ}/\bbQ), \Pic(S_{\overline\bbQ})) \longrightarrow \Br(S) / \pi^*\!\Br(\bbQ)$$
making the diagram
$$
\definemorphism{gleich}\Solid\notip\notip
\diagram
H^1(\Gal(\overline{\bbQ}/\bbQ), \Pic(S_{\overline\bbQ})) \rto^{\;\;\;\;\;\;\;\;\;\;\delta} \dto_d & \Br(S) / \pi^*\!\Br(\bbQ) \ddto^\res \\
H^2(\Gal(\overline{\bbQ}/\bbQ), \overline\bbQ(S)^*/\overline\bbQ^*) \dgleich & \\
H^2(\Gal(\overline{\bbQ}/\bbQ), \overline\bbQ(S)^*)/\pi^*\!\Br(\bbQ) \rto^{\;\;\;\;\;\;\;\;\;\;\;\;\;\inf} & \Br(\bbQ(S))/\pi^*\!\Br(\bbQ)
\enddiagram
$$
commute.
Here,~$d$
is induced by the short exact~sequence
$$0 \to \overline\bbQ(S)^*/\overline\bbQ^* \to \Div(S_{\overline\bbQ}) \to \Pic(S_{\overline\bbQ}) \to 0$$
and the other morphisms are the canonical~ones.\smallskip

\noindent
{\bf Proof.}
{\em
The~equality at the lower left corner comes from the fact~\cite[section~11.4]{Ta} that
$H^3 (\Gal(\overline{\bbQ}/\bbQ), \overline\bbQ^*) = 0$.
The~main assertion is~\cite[Lem\-ma~43.1.1]{Ma}.
}
\eop
\end{lem}

\begin{remark}
\label{fin}
The~group
$H^1(\Gal(\overline{\bbQ}/\bbQ), \Pic(S_{\overline\bbQ}))$
is always~finite. Hence, by Remark~\ref{standard}.iii), we know that only finitely many primes are relevant for the Brauer-Manin~obstruction.
\end{remark}

\section{%\boldmath The non-trivial Brauer class in the $[12,15]$-case
One double-six}

\begin{ttt}
Recall~that, on a non-singular cubic surface, there are exactly 27~lines. Their~configuration was studied intensively by the geometers of the 19th~century. They~showed that the 27~lines contain 72~{\em sixers\/} of lines which are mutually~skew. For~each sixer, there is a complementary one formed by the lines meeting exactly five lines of the~sixer. A~pair of complementary sixers is called a {\em double-six}.

The~automorphism group
$W(E_6)$
acts transitively on double-sixes. In~the blown-up model, one double-six may easily be written down~explicitly. For~this, we refer to~\cite[Remark~V.4.9.1]{Ha} from which we also adopt the notation for the~lines.
\end{ttt}

\begin{lem}
Let\/~$S$
be a non-singular cubic surface
over\/~$\bbQ$.
Suppose~that, under the operation
of\/~$\Gal(\overline\bbQ/\bbQ)$,
the 27~lines
on\/~$S$
decompose into orbits one of which is of size~15. Then,~the complementary twelve lines form a~double-six.\smallskip

\noindent
{\bf Proof.}
{\em
The~Galois
group~$\Gal(\overline\bbQ/\bbQ)$
operates via a certain
subgroup~$G \subseteq W(E_6)$.
Our~assumption implies
that~$5 | \#G$.
I.e.,~$G$~contains
the 5-Sylow~subgroup
of~$W(E_6)$.

The~operation of this is easily described in the blown-up~model. The~cyclic group
$\langle (12345) \rangle \subset S_6$
acts on the~indices. The~two lines
$E_0$
and~$G_0$
are stationary while the others form five orbits of size five~each. These~are
$E := \{E_1, \ldots, E_5\}$,
$G := \{G_1, \ldots, G_5\}$,
$F_0 := \{F_{01}, \ldots, F_{05}\}$,
$F_1 := \{F_{12}, F_{23}, F_{34}, F_{45}, F_{15}\}$,
and, finally,
$F_2 := \{F_{13}, F_{24}, F_{35}, F_{14}, F_{25}\}$.

The~intersection matrix of the five latter orbits turns out to~be
$$
\left(
\begin{array}{rrrrr}
-5 & 20 &  5 & 10 & 10 \\
20 & -5 &  5 & 10 & 10 \\
 5 &  5 & -5 & 15 & 15 \\
10 & 10 & 15 &  5 &  5 \\
10 & 10 & 15 &  5 &  5
\end{array}
\right) \, .
$$
We~assert that a size~15 orbit must be formed by
$F_0$,
$F_1$,
and~$F_2$.

Indeed,~three orbits of size five may be put together to form an orbit only if, for the corresponding divisors,
$D(D + D^\prime + D^{\prime\prime}) = D^\prime(D + D^\prime + D^{\prime\prime}) = D^{\prime\prime}(D + D^\prime + D^{\prime\prime})$.
This~excludes all combinations, except
for~$E \cup G \cup F_1$,
($E \cup G \cup F_2$,
and~$F_0 \cup F_1 \cup F_2$).
The~set
$E \cup G \cup F_1$
contains, however, only two fivers of skew lines, namely
$E$
and~$G$.
As~the lines
in~$F_1$
are missing, this is a~contradiction.
}
\eop
\end{lem}

\begin{remark}
$U_1$,
the largest subgroup
of~$W(E_6)$
stabilizing a double-six is isomorphic
to~$S_6 \times \bbZ/2\bbZ$
of order~1440.
\end{remark}

\begin{nota}
\label{f30}
Let~$S$
be a non-singular cubic~surface. Assume~that twelve of the 27~lines
on~$S$
form a double-six which
is~$\Gal(\overline\bbQ/\bbQ)$-invariant.
Choose~such a double-six. Then,~there are two kinds of tritangent~planes. We~have 15~tritangent~planes which
meet~$S$
only within the 15 complementary~lines. The~other 30~tritangent~planes contain one of the 15~lines and two from the \mbox{double-six}.

We~write
$F_{30}$
for a product over the linear~forms defining the 30~tritangent~planes and
$F_{15}$
for a product over the linear~forms defining the 15~others.
Note~that~$F_{30}/F_{15}^2 \in \bbQ(S)$.
\end{nota}

\begin{theo}
\label{eineds}
Let\/~$\pi\colon S \to \Spec \bbQ$
be a non-singular cubic surface such that the 27~lines have orbit structure\/
$[12, 15]$
under the operation
of\/~$\Gal(\overline\bbQ/\bbQ)$.

\begin{iii}
\item
Then,~$\Br(S)/\pi^*\!\Br(\bbQ) = \bbZ/2\bbZ$.
\item
For\/~$0 \neq c \in \Br(S)/\pi^*\!\Br(\bbQ)$,
a
representative\/~$\underline{c}$
of\/~$\res(c) \in \Br(\bbQ(S))/\pi^*\!\Br(\bbQ)$
is given as~follows.

Consider~the quadratic number
field\/~$\bbQ(\sqrt{D})$
splitting the double-six into two~sixers.
Then,~apply to the~class
$$(F_{30}/F_{15}^2) \in \widehat{H}^0(\Gal(\bbQ(\sqrt{D})/\bbQ), \bbQ(\sqrt{D})(S)^*) = \bbQ(S)^* / N\bbQ(\sqrt{D})(S)^*$$
the periodicity isomorphism
to\/~$H^2$
and the inflation~map.
\end{iii}\medskip\pagebreak[3]

\noindent
{\bf Proof.}
{\em
{\em First step.}
Inflation.\smallskip

\noindent
We~have the isomorphism
$\smash{\delta\colon H^1(\Gal(\overline{\bbQ}/\bbQ), \Pic(S_{\overline\bbQ})) \to \Br(S) / \pi^*\!\Br(\bbQ)}$
and will work with the group on the~left.

An~element of the
group~$\Gal(\overline{\bbQ}/\bbQ)$
may either flip the two sixers forming the twelve lines or~not. Therefore,~there is an index two subgroup stabilizing the~sixers. This~group corresponds to the quadratic number
field~$\bbQ(\sqrt{D})$.
By~Fact~\ref{blown} below, we
know~$H^1(\Gal(\overline{\bbQ}/\bbQ(\sqrt{D})), \Pic(S_{\overline\bbQ})) = 0$.
The~inflation-restriction-sequence yields~that
$$\inf\colon H^1(\Gal(\bbQ(\sqrt{D})/\bbQ), \Pic(S_{\overline\bbQ})^{\Gal(\overline\bbQ/\bbQ(\sqrt{D}))}) \longrightarrow H^1(\Gal(\overline{\bbQ}/\bbQ), \Pic(S_{\overline\bbQ}))$$
is an~isomorphism.\medskip

\noindent
{\em Second step.}
Divisors.\smallskip

\noindent
The~orbit structure of the 27~lines under the operation
of~$\Gal(\overline\bbQ/\bbQ(\sqrt{D}))$
is~$[6, 6, 15]$.
Indeed,~when going over to an index two subgroup an orbit of odd size must not~split. Denote~by
$E$,
$G$,
and~$F$
the divisors formed by summing over the first, second, and third orbit,~respectively.
$E$,
$G$,
and~$F$
clearly define elements
in~$\smash{\Pic(S_{\overline\bbQ})^{\Gal(\overline\bbQ/\bbQ(\sqrt{D}))}}$.
Write~$P$
for the subgroup generated by these three~divisors.

$P$~is
of finite index
in~$\smash{\Pic(S_{\overline\bbQ})^{\Gal(\overline\bbQ/\bbQ(\sqrt{D}))}}$.
Indeed,~every element of
$\Pic(S_{\overline\bbQ})$
is an integral linear combination of the divisors given by the 27~lines. Therefore,~every element
in~$\smash{\Pic(S_{\overline\bbQ})^{\Gal(\overline\bbQ/\bbQ(\sqrt{D}))}}$
is a
\mbox{$\bbQ$-linear}
combination
of~$E$,
$G$,
and~$F$
and the denominators are at most six or~15.

We~claim that the index
of~$P$
is a divisor of~15. In~fact, we have the relation
$5E + 5G - 4F \,\sim\, 0$.
Further,~the discriminant of the lattice spanned
by~$E$
and~$F$~is
$$\left|
\begin{array}{rr}
-6 & 30 \\
30 & 75
\end{array}
\right|
= -1350 = (-6) \!\cdot\! 15^2 \, .$$
Consequently,~$P$
is of odd index
in~$\smash{\Pic(S_{\overline\bbQ})^{\Gal(\overline\bbQ/\bbQ(\sqrt{D}))}}$.
This~implies that the natural~homomorphism
$$H^1(\Gal(\bbQ(\sqrt{D})/\bbQ), P) \longrightarrow H^1(\Gal(\bbQ(\sqrt{D})/\bbQ), \Pic(S_{\overline\bbQ})^{\Gal(\overline\bbQ/\bbQ(\sqrt{D}))})$$
is a~bijection.\medskip

\noindent
{\em Third step.}
The~fundamental~class.\smallskip

\noindent
As~$\Gal(\bbQ(\sqrt{D})/\bbQ)$
is a cyclic group of order two, we~have
$$H^2(\Gal(\bbQ(\sqrt{D})/\bbQ), \bbZ) \cong \bbZ/2\bbZ \, .$$
Write~$u$
for the non-zero~element. Then,~the periodicity isomorphism is given~by
$${} \cup u \colon \widehat{H}^{-1} (\Gal(\bbQ(\sqrt{D})/\bbQ), P) \longrightarrow H^1(\Gal(\bbQ(\sqrt{D})/\bbQ), P) \, .$$
Observe~that, for a cyclic group of order two, this isomorphism is canonical as there is no ambiguity in the choice
of~$u$.\medskip

\noindent
{\em Fourth step.}
Computing~$\widehat{H}^{-1}$.\smallskip

\noindent
We~have
$P = S/S_0$
for~$S := \bbZ E \oplus \bbZ G \oplus \bbZ F$
and~$S_0$
the group of the principal divisors contained
in~$S$.
The~relation
$\widehat{H}^{-1} (\Gal(\bbQ(\sqrt{D})/\bbQ), S) = 0$
follows immediately from the~definition. Hence,~the short exact~sequence
$$0 \to S_0 \to S \to P \to 0$$
yields
\begin{align*}
\widehat{H}^{-1} (\Gal(\bbQ(\sqrt{D})/\bbQ&), P) = \\
 & = \ker (\widehat{H}^0 (\Gal(\bbQ(\sqrt{D})/\bbQ), S_0) \to \widehat{H}^0 (\Gal(\bbQ(\sqrt{D})/\bbQ), S)) \\
 & = \ker (S_0^{\Gal(\bbQ(\sqrt{D})/\bbQ)} / N\!S_0 \to S^{\Gal(\bbQ(\sqrt{D})/\bbQ)} / N\!S) \\
 & = (S_0^{\Gal(\bbQ(\sqrt{D})/\bbQ)} \cap N\!S) / N\!S_0 \, .
\end{align*}
Here,~the norm~map acts by the rule
$N \colon aE+bG+cF \mapsto (a+b)E + (a+b)G + 2cF$.
Hence,~$N\!S = \langle E+G, 2F \rangle$.
Principal~divisors are characterized by the property that all intersection numbers are~zero. A~direct calculation~shows
$$S_0^{\Gal(\bbQ(\sqrt{D})/\bbQ)} \cap N\!S = \langle 5E + 5G - 4F \rangle \, .$$
The~generator is the norm of any divisor of the
form~$aE + (5-a)G - 2F$.
None~of these is principal. Indeed,~the intersection number
with~$E$
is equal to
$-6a + 30(5-a) - 60 = -36a + 90$
and this terms does not vanish
for~$a \in \bbZ$.
Assertion~i) is~proven.\medskip

\noindent
{\em Fifth step.}
The~representative.\smallskip

\noindent
We~actually constructed a non-zero
element~$c^\prime \in H^1(\Gal(\overline{\bbQ}/\bbQ), \Pic(S_{\overline\bbQ}))$.
The~calculations given above show that
$$d(c^\prime) = (F_{30}/F_{15}^2) \cup u \, .$$
Indeed,~it is easy to see
that~$\div (F_{30}/F_{15}^2) = 5E + 5G - 4F$.
The~assertion now follows from the commutative diagram given in~Lemma~\ref{75}.
}
\eop
\end{theo}

\begin{fac}
\label{blown}
Let\/~$S$
be a non-singular cubic surface over a
field~$K$
obtained by blowing
up\/~$\smash{\bP^2_K}$
in six\/
\mbox{$\overline{K}$-rational}
points which form a Galois invariant~set.\smallskip

\noindent
Then,~$\smash{H^1(\Gal(\overline{K}/K), \Pic(S_{\overline{K}})) = 0}$.\medskip

\noindent
{\bf Proof.}
{\em
This~is, of course, a particular case of the very general~\cite[Theorem~29.1]{Ma}. Let~us give an elementary proof~here.

According~to Shapiro's~lemma, we may replace
$\Gal(\overline{K}/K)$
by a finite
quotient~$G$.
We~have
$\smash{\Pic(S_{\overline{K}}) = \bbZ H \oplus \bbZ E_1 \oplus \ldots \oplus \bbZ E_6}$
for~$H$
the hyperplane section
of~$\bP^2$
and
$E_1, \ldots, E_6$
the exceptional~divisors. Therefore,~as
a~\mbox{$G$-module,}
$$\Pic(S_{\overline{K}}) = \bbZ \oplus \bbZ[G/H_1] \oplus \ldots \oplus \bbZ[G/H_l]$$
for
$l$
the number of Galois~orbits and certain
subgroups~$H_1, \ldots, H_l$.
Clearly,~we have
$H^1 (G, \bbZ) = \Hom (G, \bbZ) = 0$.

On~the other hand, for any
subgroup~$H$,
the
\mbox{$G$-module~$\bbZ[G/H]$}
is equipped with a non-degenerate~pairing.
Hence,~$\bbZ[G/H] \cong \Hom(\bbZ[G/H], \bbZ)$
and
\begin{eqnarray*}
H^1 (G, \bbZ[G/H]) & \cong & H^1 (G, \Hom(\bbZ[G/H], \bbZ)) \\
 & \cong & \smash{\widehat{H}^0 (G, \Hom(\bbZ[G/H], \bbQ/\bbZ))} \\
 & \cong & \smash{\Hom\!\big(\widehat{H}^{-1} (G, \bbZ[G/H]), \bbQ/\bbZ\big)}
\end{eqnarray*}
by~the duality theorem for cohomology of finite groups~\cite[Chap.~XII, Corollary~6.5]{CE}.
Finally,~$\smash{\widehat{H}^{-1} (G, \bbZ[G/H])}$
vanishes as is seen immediately from the~definition.
}
\eop
\end{fac}

\section{Triples of azygetic double-sixes%---\\The~$[6, 6, 6, 9]$-case
}

\begin{defi}
\label{class}
Let~$\pi\colon S \to \Spec\bbQ$
be a smooth cubic surface and
$\calD$~be
a Galois invariant double-six. This~induces a group~homomorphism
$\Gal(\overline\bbQ/\bbQ) \to U_1$,
given by the operation on the 27~lines.

\begin{iii}
\item
Then,~the image of the non-zero element under the natural~homomorphism
$$\bbZ/2\bbZ \cong H^1(U_1, \Pic(S_{\overline\bbQ})) \longrightarrow H^1(\Gal(\overline\bbQ/\bbQ), \Pic(S_{\overline\bbQ})) \stackrel{\delta}{\longrightarrow} \Br(S)/\pi^*\!\Br(\bbQ)$$
is called {\em the Brauer class associated with\/} the
double-six~$\calD$.
We~denote it
by~$\cl(\calD)$.
\item
This~defines a map
$$\cl\colon \Phi_S^{\Gal(\overline\bbQ/\bbQ)} \longrightarrow \Br(S)/\pi^*\!\Br(\bbQ)$$
from the set of all Galois invariant double-sixes which we will call the {\em class map}.
\end{iii}
\end{defi}

\begin{ttt}
\label{syzy}
Two~double-sixes may have either four or six lines in~common. In~the former case, the two are called syzygetic, in the latter~azygetic. A~pair of azygetic double-sixes is built as~follows.
$$
\left(
\begin{array}{cccccc}
E_0 & E_1 & E_2 & E_3 & E_4 & E_5 \\
G_0 & G_1 & G_2 & G_3 & G_4 & G_5 
\end{array}
\right),
\quad
\left(
\begin{array}{cccccc}
E_0    & E_1    & E_2    & F_{45} & F_{35} & F_{34} \\
F_{12} & F_{02} & F_{01} & G_3    & G_4    & G_5 
\end{array}
\right).
$$
The~twelve~lines which appear only once form a third \mbox{double-six}
$$
\left(
\begin{array}{cccccc}
F_{12} & F_{02} & F_{01} & E_3    & E_4    & E_5    \\
G_0    & G_1    & G_2    & F_{45} & F_{35} & F_{34} 
\end{array}
\right)
$$
azygetic to the other~two.
\end{ttt}

\begin{ttt}
The~largest
group~$U_3$
stabilizing a triple of azygetic double-sixes is isomorphic
to~$(S_3 \times S_3) \ltimes \bbZ/2\bbZ$
of order~72. The~induced orbit structure
is~$[6, 6, 6, 9]$.
The~orbits themselves are, in the notation above,
$\{E_0, E_1, E_2, G_3, G_4, G_5\}$,
$\{G_0, G_1, G_2, E_3, E_4, E_5\}$,
$\{F_{01}, F_{02}, F_{12}, F_{34}, F_{35}, F_{45}\}$,
and
$\{F_{03}, F_{04}, F_{05}, F_{13}, F_{14}, F_{15}, F_{23}, F_{24}, F_{25}\}$.

The~quadratic extension
$\bbQ(\sqrt{D})$
splitting one of the three double-sixes into two sixers automatically splits the others,~too. The~operation
of~$\Gal(\overline\bbQ/\bbQ(\sqrt{D}))$
yields the orbits
$\{E_0, E_1, E_2\}$,
$\{E_3, E_4, E_5\}$,
$\{G_0, G_1, G_2\}$,
$\{G_3, G_4, G_5\}$,
$\{F_{01}, F_{02}, F_{12}\}$,
$\{F_{34}, F_{35}, F_{45}\}$,
and
$\{F_{03}, F_{04}, F_{05}, F_{13}, F_{14}, F_{15}, F_{23}, F_{24}, F_{25}\}$.
\end{ttt}

\begin{theo}
\label{3syzy}
Let\/~$\pi\colon S \to \Spec \bbQ$
be a non-singular cubic~surface. Assume that\/
$\Gal(\overline\bbQ/\bbQ)$
stabilizes a
triple~$\{\calD_1, \calD_2, \calD_3\}$
of azygetic double-sixes and that the 27~lines have orbit structure\/~$[6, 6, 6, 9]$.

\begin{iii}
\item
Then,~$\Br(S)/\pi^*\!\Br(\bbQ) = \bbZ/2\bbZ \times \bbZ/2\bbZ$.
\item
The~three non-zero elements are\/
$\cl(\calD_1)$,
$\cl(\calD_2)$,
and\/~$\cl(\calD_3)$.
\end{iii}\smallskip

\noindent
{\bf Proof.}
{\em
{\em First step.}
Inflation.\smallskip

\noindent
We~will follow the same strategy as in the proof of Theorem~\ref{eineds}. In~particular, we will work with the
group~$H^1(\Gal(\bbQ(\sqrt{D})/\bbQ), \Pic(S_{\overline\bbQ})^{\Gal(\overline\bbQ/\bbQ(\sqrt{D}))})$.\medskip

\noindent
{\em Second step.}
Divisors.\smallskip

\noindent
The~orbit structure of the 27~lines under the operation
of~$\Gal(\overline\bbQ/\bbQ(\sqrt{D}))$
is
$[3, 3, 3, 3, 3, 3, 9]$.
For~the lines and double-sixes, we use the notation introduced in~\ref{syzy}. Further,~denote~by
$$E^{(1)} ,E^{(2)}, G^{(1)}, G^{(2)}, F^{(1)}, F^{(2)}, F^{(3)}$$
the divisors formed by summing over the~orbits. These~clearly define elements in
$\smash{\Pic(S_{\overline\bbQ})^{\Gal(\overline\bbQ/\bbQ(\sqrt{D}))}}$.
We~write~$P$
for the subgroup generated by these seven~divisors.

Every~element of
$\Pic(S_{\overline\bbQ})$
is an integral linear combination of the divisors given by the 27~lines. Therefore,~every element
in~$\smash{\Pic(S_{\overline\bbQ})^{\Gal(\overline\bbQ/\bbQ(\sqrt{D}))}}$
is a
\mbox{$\bbQ$-linear}
combination
of~$E^{(1)}$,
$E^{(2)}$,
$G^{(1)}$,
$G^{(2)}$,
$F^{(1)}$,
$F^{(2)}$,
and~$F^{(3)}$
and the denominators are divisors of~nine.
Consequently,~$P$
is of odd index
in~$\smash{\Pic(S_{\overline\bbQ})^{\Gal(\overline\bbQ/\bbQ(\sqrt{D}))}}$.
This~implies that the natural~homomorphism
$$H^1(\Gal(\bbQ(\sqrt{D})/\bbQ), P) \longrightarrow H^1(\Gal(\bbQ(\sqrt{D})/\bbQ), \Pic(S_{\overline\bbQ})^{\Gal(\overline\bbQ/\bbQ(\sqrt{D}))})$$
is a~bijection.\medskip

\noindent
{\em Third step.}
The~fundamental~class.\smallskip

\noindent
Again,~we
write~$u$
for the non-zero~element
in~$H^2(\Gal(\bbQ(\sqrt{D})/\bbQ), \bbZ) \cong \bbZ/2\bbZ$.
Then,~the periodicity isomorphism is given~by
$${} \cup u \colon \widehat{H}^{-1} (\Gal(\bbQ(\sqrt{D})/\bbQ), P) \longrightarrow H^1(\Gal(\bbQ(\sqrt{D})/\bbQ), P) \, .$$

\noindent
{\em Fourth step.}
Computing~$\widehat{H}^{-1}$.\smallskip

\noindent
We~have
$P = S/S_0$
for~$S := \bbZ E^{(1)} \oplus \bbZ E^{(2)} \oplus \bbZ G^{(1)} \oplus \bbZ G^{(2)} \oplus \bbZ F^{(1)} \oplus \bbZ F^{(2)} \oplus \bbZ F^{(3)}$
and~$S_0$
the group of the principal divisors contained
in~$S$.
To~simplify formulas, we will use the notation
$D^{(1)} := E^{(1)}\!+\!G^{(2)}$,
$D^{(2)} := E^{(2)}\!+\!G^{(1)}$,
and~$D^{(3)} := F^{(1)}\!+\!F^{(2)}$.

The~relation
$\widehat{H}^{-1} (\Gal(\bbQ(\sqrt{D})/\bbQ), S) = 0$
follows immediately from the~definition. Hence,~the short exact~sequence
$0 \to S_0 \to S \to P \to 0$
yields, as~above,
$$\widehat{H}^{-1} (\Gal(\bbQ(\sqrt{D})/\bbQ), P) = (S_0^{\Gal(\bbQ(\sqrt{D})/\bbQ)} \cap N\!S) / N\!S_0 \, .$$
Here,~the norm~map acts by the rule
\begin{align*}
& N \colon a_1E^{(1)} + a_2E^{(2)} + b_1G^{(1)} + b_2G^{(2)} + c_1F^{(1)} + c_2F^{(2)} + c_3F^{(3)} \\
& \hspace{20mm} \mapsto (a_1+b_2) D^{(1)} + (a_2+b_1) D^{(2)} + (c_1+c_2) D^{(3)} + 2c_3 F^{(3)} \, .
\end{align*}
Hence,~$N\!S = \langle D^{(1)}, D^{(2)}, D^{(3)}, 2F^{(3)}\rangle$.
Principal~divisors are characterized by the property that all intersection numbers are~zero. A~direct calculation~shows
\begin{eqnarray*}
\!\!\! &   & \smash{S_0^{\Gal(\bbQ(\sqrt{D})/\bbQ)} \cap N\!S} \\
\!\!\! & = & \smash{\{ d_1 D^{(1)} + d_2 D^{(2)} + d_3 D^{(3)} + 2 e F^{(3)} \mid d_1 + d_2 + d_3 + 3 e = 0} \} \\
\!\!\! & = & \langle D^{(1)} - D^{(2)}, D^{(1)} - D^{(3)}, D^{(1)} + D^{(2)} + D^{(3)} - 2F^{(3)} \rangle \, .
\end{eqnarray*}
It~is easy to see that
$E^{(1)} + G^{(1)} + F^{(1)} - F^{(3)}$
is a principal~divisor.
Hence,
$$D^{(1)} + D^{(2)} + D^{(3)} - 2F^{(3)} \in N\!S_0 \, .$$
Further,~$N\!S_0$
contains all principal divisors which are divisible
by~$2$.
As~it turns out that these two sorts of elements generate the whole
of~$N\!S_0$,
assertion~i)~follows.\medskip

\noindent
{\em Fifth step.}
The~representatives.\smallskip

\noindent
We~actually constructed non-zero\vspace{0.2mm}
elements~$\smash{c_1, c_2, c_3 \in H^1(\Gal(\overline{\bbQ}/\bbQ), \Pic(S_{\overline\bbQ}))}$,
represented~by
$\smash{D^{(1)} - D^{(2)}, D^{(1)} - D^{(3)}, D^{(2)} - D^{(3)} \in S_0^{\Gal(\bbQ(\sqrt{D})/\bbQ)} \cap N\!S}$.
The~first representative is equivalent to
$$3(D^{(1)} \!-\! D^{(2)}) + 2(D^{(1)} \!+\! D^{(2)} \!+\! D^{(3)}\! -\! 2F^{(3)}) + 6(D^{(2)} \!-\! D^{(3)}) = \div (F_{30}/F_{15}^2) \, .$$
Hence,~$d(c_1) = (F_{30}/F_{15}^2) \cup u$.
The~corresponding Brauer class
is~$\cl\!\big(({E_0 \ldots E_5 \atop G_0 \ldots G_5})\big)$.
For~the two other classes, the situation is~analogous.
}
\eop
\end{theo}

\begin{remark}
In~the
$[6, 6, 6, 9]$-case,
the 45~tritangent planes decompose into five~orbits.\smallskip

\begin{iii}
\item[$\bullet$ ]
$[E_i, G_j, F_{ij}]$
for~$i,j \in \{0,1,2\}$
or~$i,j \in \{3,4,5\}$,
$i\neq j$.\quad
(twelve~planes)
\item[$\bullet$ ]
$[E_i, G_j, F_{ij}]$
for~$i \in \{0,1,2\}$
and~$j \in \{3,4,5\}$.\quad
(nine~planes)
\item[$\bullet$ ]
$[E_i, G_j, F_{ij}]$
for~$i \in \{3,4,5\}$
and~$j \in \{0,1,2\}$.\quad
(nine~planes)
\item[$\bullet$ ]
$[F_{i_0 i_1}, F_{i_2 i_3}, F_{i_4 i_5}]$
for~$\{i_0, i_1, i_2, i_3, i_4, i_5\} = \{0, 1, 2, 3, 4, 5\}$,\hfill\break
\begin{minipage}{5.5cm}
~
\end{minipage}
$i_0,i_1 \in \{0,1,2\}$,
and
$i_2,i_3 \in \{3,4,5\}$.\quad
(nine~planes)
\item[$\bullet$ ]
$[F_{i_0 i_1}, F_{i_2 i_3}, F_{i_4 i_5}]$
for~$\{i_0, i_2, i_4\} = \{0,1,2\}$
and~$\{i_1, i_3, i_5\} = \{3,4,5\}$.\quad
(six~planes)
\end{iii}\smallskip

\noindent
The~three forms of
type~$F_{30}$
are obtained by multiplying the linear forms defining the orbit of size twelve together with those for two of the orbits of size~nine. Actually,~the size twelve orbit is~irrelevant. Up~to a scalar factor, it is the square of a sextic~form.
\end{remark}

\begin{remark}
Triples~of azygetic double-sixes have been studied by the classical algebraic~geometers. See,~for example, \cite[\S6]{Ko}. A~result from the 19th~century states that there are exactly 120 triples of azygetic double-sixes on a smooth cubic~surface. Actually,~the automorphism group
$W(E_6)$
acts transitively on~them.
%The~stabilizer is the group of order~72 given~above.
\end{remark}

\section{The general case of a Galois group stabilizing\\ a double-six}

\begin{ttt}
To~explicitly compute
$H^1 (G, \Pic(S_{\overline\bbQ}))$
as an abstract abelian group, one may use Manin's formula~\cite[Proposition~31.3]{Ma}.
This~means the~following.

$\smash{\Pic(S_{\overline\bbQ})}$~is
generated by the 27~lines. The~group of all permutations of the 27~lines respecting the intersection pairing is isomorphic to the Weyl
group~$W(E_6)$
of
order~$51\,840$.
The~group~$G$
operates on the 27~lines via a finite
quotient~$G/H$
which is isomorphic to a subgroup
of~$W(E_6)$.
Then,
$$H^1 (G, \Pic(S_{\overline\bbQ})) \cong \Hom(N\!D \cap D_0 / N\!D_0, \bbQ/\bbZ) \, .$$
Here,~$D$
is the free abelian group generated by the 27~lines and
$D_0$
is the subgroup of all principal~divisors.
$N \colon D \to D$
denotes the norm map as a
\mbox{$G/H$-module}.
\end{ttt}

\begin{ttt}
Using~Manin's formula, we computed
$H^1 (G, \Pic(S_{\overline\bbQ}))$
for each of the 350~conjugacy classes of subgroups
of~$W(E_6)$.
The~computations in~{\tt GAP} took approximately 28~seconds of CPU~time. Thereby,~we recovered the following result of Sir~Peter~Swinnerton-Dyer~\cite{SD}. (See also~P.~K.~Corn~\cite{Co}.)
\end{ttt}

\begin{theoo}[{\rm Swinnerton-Dyer}{}]
Let\/~$S$
be a non-singular cubic surface
over\/~$\bbQ$.
Then,~$H^1 (\Gal(\overline\bbQ/\bbQ), \Pic(S_{\overline\bbQ}))$
may take only five values,
$0$,
$\bbZ/2\bbZ$,
$\bbZ/2\bbZ \times \bbZ/2\bbZ$,
$\bbZ/3\bbZ$,
and~$\bbZ/3\bbZ \times \bbZ/3\bbZ$.
\eop
\end{theoo}

\begin{remark}
\label{GAP}
One~has
$H^1 (G, \Pic(S_{\overline\bbQ})) = 0$
in 257 of the 350~cases.
\end{remark}

\begin{ttt}
\label{exp}
More~importantly, we make the following observation.\smallskip

\noindent
{\bf Proposition.}
{\em
Let\/~$S$
be a non-singular cubic surface
over\/~$\bbQ$.

\begin{iii}
\item
If\/~$\smash{H^1 (\Gal(\overline\bbQ/\bbQ), \Pic(S_{\overline\bbQ})) = \bbZ/2\bbZ}$
then,~on\/~$S$,
there is a Galois invariant double-six.
\item
If\/~$\smash{H^1 (\Gal(\overline\bbQ/\bbQ), \Pic(S_{\overline\bbQ})) = \bbZ/2\bbZ \times \bbZ/2\bbZ}$
then,~on\/~$S$,
there is a triple of azy\-getic double-sixes stabilized
by~$\Gal(\overline\bbQ/\bbQ)$.
\end{iii}
\smallskip

\noindent
{\bf Proof.}
{\em
This~is seen by a case-by-case study using~{\tt GAP}.
}}
\eop
\end{ttt}

\begin{rems}
\label{12oder4}
\begin{iii}
\item
On~the other hand, if there is a Galois invariant double-six
on~$S$
then
$H^1 (\Gal(\overline\bbQ/\bbQ), \Pic(S_{\overline\bbQ}))$
is either
$0$,
or~$\bbZ/2\bbZ$
or~$\bbZ/2\bbZ \times \bbZ/2\bbZ$.
\item
Proposition~\ref{exp} immediately provokes the question whether the cohomology classes are always ``the same'' as in the
$[12, 15]$-
and~$[6, 6, 6, 9]$-cases.
I.e.,~of the type
$\cl(\calD)$
for certain Galois invariant double-sixes. Somewhat~surprisingly, this is indeed the~case.
\end{iii}
\end{rems}

\begin{lem}
\label{criteria}
Let\/~$\calS$
be a non-singular cubic surface over an algebraically closed field,
$H$~a
group of automorphisms of the configuration of the 27~lines, and\/
$H^\prime \subseteq H$
any~subgroup. Each~of the criteria below is sufficient~for
$$\res\colon H^1 (H, \Pic(\calS)) \to H^1 (H^\prime, \Pic(\calS))$$
being an~injection.\smallskip

\begin{iii}
\item
$H$
and\/
$H^\prime$
generate the same orbit~structure.
\item
$H^\prime$~is
of odd index
in\/~$H$
and\/
$H^1 (H, \Pic(\calS))$~is
a\/~\mbox{$2$-group.}
\item
$H^\prime$~is
a normal subgroup
in\/~$H$
and\/~$\rk \Pic(\calS)^H = \rk \Pic(\calS)^{H^\prime}$.
\end{iii}\smallskip

\noindent
{\bf Proof.}
{\em
i)
This~follows immediately from the formula of~Manin~\cite[Proposition~31.3]{Ma}.\smallskip

\noindent
ii)
Here,~$\cores \circ \res \colon H^1 (H, \Pic(\calS)) \to H^1 (H, \Pic(\calS))$
is the multiplication by an odd number, hence the~identity.\smallskip

\noindent
iii)
The~assumption ensures that
$H/H^\prime$~operates
trivially
on~$\Pic(\calS)^{H^\prime}$.
Hence,
$H^1 \!\big( H/H^\prime, \Pic(\calS)^{H^\prime} \big) = 0$.
The~inflation-restriction~sequence
$$0 \to H^1 \!\big( H/H^\prime, \Pic(\calS)^{H^\prime} \big) \to H^1(H, \Pic(\calS)) \to H^1(H^\prime, \Pic(\calS))$$
yields the~assertion.
}
\eop
\end{lem}

\begin{prop}
Let\/~$\calS$
be a non-singular cubic surface over an algebraically closed field,
$U_1$~the
group of automorphisms of the configuration of the 27~lines stabilizing a double-six and\/
$U_3$~the
group stabilizing a triple of azygetic double-sixes.

\begin{abc}
\item
Let\/~$H \subseteq U_1$
be such
that\/~$H^1 (H, \Pic(\calS)) = \bbZ/2\bbZ$.
Then,~the restriction
$$\res\colon H^1 (U_1, \Pic(\calS)) \longrightarrow H^1 (H, \Pic(\calS))$$
is a~bijection.
\item
Let\/~$H \subseteq U_3$
be a subgroup such
that\/~$H^1 (H, \Pic(\calS)) = \bbZ/2\bbZ \times \bbZ/2\bbZ$.
Then,~the restriction
$$\res\colon H^1 (U_3, \Pic(\calS)) \longrightarrow H^1 (H, \Pic(\calS))$$
is a~bijection.
\end{abc}\smallskip

\noindent
{\bf Proof.}
{\em
The~proof has a computer~part. We~use the machine to verify that the criteria formulated in Lemma~\ref{criteria} suffice to establish the result in all~cases.\smallskip

\noindent
b) Here,~the subgroup
$(A_3 \times A_3) \ltimes \bbZ/2\bbZ$
of order~18, as well as the two intermediate groups of order~36 produce the same orbit
structure~$[6, 6, 6, 9]$.
It~turns out that every
subgroup~$H$
which leads
to~$\bbZ/2\bbZ \times \bbZ/2\bbZ$
is a subgroup of odd
\mbox{($1$,
$3$,
or~$9$)}
index in one of~those.\smallskip

\noindent
a)
By~Lemma~\ref{criteria}.ii), we may test this on the
\mbox{$2$-Sylow~subgroups}
of~$H$
and~$U_1$.
$U_1^{(2)}$~is
a group of
order~$32$
such that the Picard rank is equal to~two. It~turns out that,
for~$2$-groups~$H^\prime$
such that 
$H^1 (H^\prime, \Pic(\calS)) = \bbZ/2\bbZ$,
the Picard rank may be only two or~three.

There~is a maximal
\mbox{$2$-group}
such that the Picard rank is~three. This~is a group of order~16 generating the orbit
structure~$[1, 1, 1, 4, 4, 4, 4, 4, 4]$.
To~prove the assertion for this group, one first observes that it is of index three in a group of order~48 with orbit
structure~$[3, 12, 12]$.
This~group, in turn, is of index two in the maximal group with that orbit~structure. That~one, being of order~96, is the maximal group stabilizing a double-six and a tritangent plane containing three complementary~lines. It~is of index~15
in~$U_1$.
}
\eop
\end{prop}

\begin{coro}
\label{inj}
Let\/~$H^\prime \subseteq H \subseteq U_1$
be~arbitrary. Then,~for the restriction map\/
$\res\colon H^1 (H, \Pic(\calS)) \longrightarrow H^1 (H^\prime, \Pic(\calS))$,
there are the following~limitations.

\begin{iii}
\item
If\/~$H^1 (H, \Pic(\calS)) = 0$
then\/~$H^1 (H^\prime, \Pic(\calS)) = 0$.
\item
If\/~$H^1 (H, \Pic(\calS)) \cong \bbZ/2\bbZ$
and\/~$H^1 (H^\prime, \Pic(\calS)) \neq 0$
then\/
$\res$
is an~injection.
\item
If\/~$H^1 (H, \Pic(\calS)) \cong \bbZ/2\bbZ \times \bbZ/2\bbZ$
then\/
$H^1 (H^\prime, \Pic(\calS)) \cong \bbZ/2\bbZ \times \bbZ/2\bbZ$
or\/~$0$.
In~the former case,
$\res$
is a~bijection. The~latter is possible only when\/
$H^\prime$
stabilizes a~sixer.
\end{iii}\smallskip

\noindent
{\bf Proof.}
{\em
We~know from Remark~\ref{12oder4}.i) that both~groups may be only
$0$,
$\bbZ/2\bbZ$,
or
$\bbZ/2\bbZ \times \bbZ/2\bbZ$.\smallskip

\noindent
i)
If~$H^1 (H^\prime, \Pic(\calS))$
were isomorphic
to~$\bbZ/2\bbZ$
or
$\bbZ/2\bbZ \times \bbZ/2\bbZ$
then the restriction from
$U_1$,
respectively~$U_3$,
to~$H^\prime$
would be the zero~map.\smallskip

\noindent
ii) is immediate from the computations~above.\smallskip

\noindent
iii)
If~$H^1 (H^\prime, \Pic(\calS)) \cong \bbZ/2\bbZ$
then, by composition, we could produce the zero map
on~$\bbZ/2\bbZ$.
The~final assertion is an experimental~observation.
}
\eop
\end{coro}

\begin{remo}[{\rm Pairs of syzygetic double-sixes}{}]
$U_2$,
the largest group stabilizing two syzygetic double-sixes is of
order~$96$.
In~view of Proposition~\ref{exp}.ii), this ensures that
$H^1 (U_2, \Pic(\calS)) \cong \bbZ/2\bbZ$
or~$0$.
Actually,~it is isomorphic
to~$\bbZ/2\bbZ$.
Corollary~\ref{inj}.ii) implies that the Brauer classes associated with the two double-sixes~coincide. Both~are equal to the non-zero~element.

Actually,~the~group~$U_2$
leads to an orbit
structure~$[1,4,6,8,8]$.
It~is easy to compute
$H^1 (U_2, \Pic(\calS))$
directly using the same methods as in the proof of Theorem~\ref{3syzy}.
\end{remo}\pagebreak[3]

\begin{theo}
Let\/~$\pi\colon S \to \Spec\bbQ$
be an arbitrary smooth cubic~surface. Then,~consider the class map
$$\cl\colon \Phi_S^{\Gal(\overline\bbQ/\bbQ)} \to \Br(S)/\pi^*\!\Br(\bbQ) \, ,$$
introduced in~Definition~\ref{class}.ii), from the set of all Galois invariant double-sixes.\smallskip

\begin{abc}
\item
$\cl$~has
the properties~below.
\begin{iii}
\item
If\/~$\calD_1, \calD_2$
are syzygetic double-sixes then\/
$\cl(\calD_1) = \cl(\calD_2)$.
\item
If\/~$\calD_1, \calD_2$
are azygetic then\/
$\cl(\calD_1) + \cl(\calD_2) + \cl(\calD_3) = 0$
for\/
$\calD_3$
the third double-six of the corresponding~triple.\smallskip
\end{iii}

\item
If\/~$\Br(S)/\pi^*\!\Br(\bbQ) \neq 0$
then
\begin{iii}
\item
\vskip-\smallskipamount
$\cl(\calD) \neq 0$
for every Galois invariant double-six. Further,
$\cl\!\big(\Phi_S^{\Gal(\overline\bbQ/\bbQ)}\big)$
contains exactly the elements of order~two.
\item
Two~double-sixes\/
$\calD_1, \calD_2$
are syzygetic if and only if\/
$\cl(\calD_1) = \cl(\calD_2)$
and azygetic if and only if\/
$\cl(\calD_1) \neq \cl(\calD_2)$.
\eop
\end{iii}
\end{abc}
\end{theo}

\section{Explicit Brauer-Manin~obstruction
%in the $[12,15]$-case
}
\label{fuenf}

\begin{ttt}
Let~$S$
be a non-singular cubic surface with a Galois invariant \mbox{double-six}~$\calD$.
This~determines a class
$c := \cl(\calD) \in \Br(S) / \pi^*\!\Br(\bbQ)$.
Choose~a representative
$\underline{c} \in \Br(S)$
and the corresponding rational
function~$F_{30}/F_{15}^2 \in \bbQ(S)$.
Finally,~let
$\bbQ(\sqrt{D})$
be the quadratic extension splitting the double-six into two~sixers.
\end{ttt}

\begin{fac}
The~quaternion algebra
over\/~$\bbQ(S)$
corresponding
to\/~$\underline{c}$~is
$$Q := \bbQ(S)\{X, Y\} / (XY + Y\!X, X^2 - D, Y^2 - F_{30}/F_{15}^2) \, .$$
\end{fac}

\begin{remark}
It~is well known that a class
in~$\Br(S)$
is uniquely determined by its restriction
to~$\Br(\bbQ(S))$.
The~corresponding quaternion algebra over the whole
of~$S$
may be described as~follows.\smallskip

Let~$x \in S$.
We~know that
$\div(F_{30}/F_{15}^2)$
is the norm of a divisor
on~$\smash{S_{\bbQ(\sqrt{D})}}$.
That~one~is necessarily locally~principal. I.e.,~we have a rational function
$\smash{f_x = a_x + b_x \sqrt{D}}$
such that
$\div(N\!f_x) = \div(F_{30}/F_{15}^2)$
on a Zariski neighbourhood
of~$x$.
Over~the maximal such
neighbourhood~$U_x$,
we define a quaternion algebra~by
$$\smash{\textstyle Q_x := \calO_{U_x}\{X, Y_x\} / (XY_x + Y_xX, X^2 - D, Y_x^2 - \frac{F_{30}}{F_{15}^2 N\!f_x}) \, .}$$
In~particular, in a
neighbourhood~$U_\eta$
of the generic point, we have the quaternion algebra
$Q_\eta := \calO_{U_\eta}\{X, Y\} / (XY + Y\!X, X^2 - D, Y^2 - F_{30} / F_{15}^2)$.

Over~$U_\eta \cap U_x$,
there is the isomorphism
$\iota_{\eta,x} \colon Q_\eta |_{U_\eta \cap U_x} \to Q_x |_{U_\eta \cap U_x}$,
given~by
$$X \mapsto X, \quad Y \mapsto (a_x + b_x X) Y_x \, .$$
For~two
points~$x, y \in S$,
the isomorphism
$\iota_{\eta,y} \!\circ\! \iota_{\eta,x}^{-1} \colon Q_x |_{U_\eta \cap U_x \cap U_y} \to Q_y |_{U_\eta \cap U_x \cap U_y}$
extends
to~$U_x \cap U_y$.

Hence,~the quaternion algebras
$Q_x$
may be glued together along these~isomorphisms. This~yields a quaternion algebra
$\calQ$
over~$S$.
\end{remark}

\begin{coro}
\label{Brexp}
Let\/~$\pi\colon S \to \Spec \bbQ$
be a non-singular cubic surface with a Galois invariant
double-six\/~$\calD$.
Further,~let\/~$p$
be a prime number and\/
$\underline{c} \in \Br(S)$
a representative of the
class\/~$\cl(\calD) \in \Br(S)/\pi^*\!\Br(\bbQ)$.\smallskip

\noindent
Then,~the local evaluation~map\/
$\ev_p (\underline{c}, \,.\; ) \colon S(\bbQ_p) \to \bbQ/\bbZ$
is given as~follows.

\begin{iii}
\item
Let\/~$x \in S(\bbQ_p)$.
Choose~a rational
function\/~$f_x$
such that\/
$\div(N\!f_x) = \div(F_{30}/F_{15}^2)$.
Then,
$$\ev_p (\underline{c}, x) = \left\{
\begin{array}{ll}
0 & {\it ~if~}
\frac{F_{30}}{F_{15}^2 N\!f_x}(x) \in \bbQ_p^*
{\it ~is~in~the~image~of\/~}
N\colon \bbQ_p(\sqrt{D}) \longrightarrow \bbQ_p \, , \\
\frac12 & {\it ~otherwise.}
\end{array}
\right.$$
Here,~$F_{30}/F_{15}^2 \in \bbQ(S)$
is the rational function corresponding to the
representative\/~$\underline{c}$.
$\bbQ(\sqrt{D})$~is
the quadratic field splitting the double-six into two~sixers.

\item
If\/~$x$
is not contained in any of the 27~lines then\/
$f_x \equiv 1$
is~allowed.
\end{iii}\smallskip

\noindent
{\bf Proof.}
{\em
Assertion~i) immediately follows from the~above. For~ii), recall that
$\div(F_{30}/F_{15}^2)$
is a linear combination of the 27~lines.
}
\eop
\end{coro}

\begin{ttt}
As~already noticed in Remark~\ref{fin}, the local evaluation map carries information only at finitely many~primes. To~exclude a particular prime, the following elementary criteria are highly~practical.
\end{ttt}

\begin{lemo}[{\rm The local $H^1$-criterion}{}]
Let\/~$S$
be a non-singular cubic surface
over\/~$\bbQ$.
Suppose~that,~for the decomposition~group
$G_p \cong \Gal(\overline\bbQ_p/\bbQ_p)$
at a prime
number\/~$p$,
$$H^1 (G_p, \Pic(S_{\overline\bbQ})) = 0 \, .$$
Then,~for every
$\alpha \in \Br(S)$,
the~value
of\/~$\ev_p (\alpha, x)$
is independent
of\/~$x \in S(\bbQ_p)$.\smallskip

\noindent
{\bf Proof.}
{\em
The~local evaluation
map~$\ev_p$
factors via
$\Br(S \times_{\Spec \bbQ} \Spec \bbQ_p)$.
By~\cite[Lemma~43.1.1]{Ma}, we have~that
$$\smash{\Br(S \times_{\Spec \bbQ} \Spec \bbQ_p) / \Br(\bbQ_p) \cong H^1 (\Gal(\overline\bbQ_p/\bbQ_p), \Pic(S_{\overline\bbQ_p}))} \, .$$
Together~with the assumption, this~yields
$\Br(S \times_{\Spec \bbQ} \Spec \bbQ_p) = \Br(\bbQ_p) = \bbQ/\bbZ$.
The~assertion~follows.
}
\eop
\end{lemo}

\begin{remark}
Recall~from Remark~\ref{GAP} that we have 
$H^1 (G, \Pic(S_{\overline\bbQ})) = 0$
for 257 of the 350 possible conjugacy classes of~subgroups.
\end{remark}

\begin{coro}
Let\/~$\pi\colon S \to \Spec \bbQ$
be a non-singular cubic surface with a Galois invariant
double-six\/~$\calD$.
Further,~let\/
$\underline{c} \in \Br(S)$
be a representative of the
class\/~$\cl(\calD) \in \Br(S)/\pi^*\!\Br(\bbQ)$.\smallskip

\noindent
If~a prime
number\/~$p$
splits in the quadratic number
field\/~$\bbQ(\sqrt{D})$
splitting the two sixers then the~value
of\/~$\ev_p (\underline{c}, x)$
is independent
of\/~$x \in S(\bbQ_p)$.\smallskip

\noindent
{\bf Proof.}
{\em
This~criterion is, of~course, an immediate consequence of Corollary~\ref{Brexp}. In~view of Fact~\ref{blown}, it is also a particular case of the local
$H^1$-criterion.
}
\eop
\end{coro}

\begin{prop}
\label{unverzw}
Let\/~$\pi\colon S \to \Spec \bbQ$
be a non-singular cubic surface with a Galois invariant
double-six\/~$\calD$.
Further,~let\/
$\underline{c} \in \Br(S)$
be a representative of the class\/
$\cl(\calD) \in \Br(S)/\pi^*\!\Br(\bbQ)$.\smallskip\pagebreak[3]

\noindent
Then,~for a prime
number\/~$p$
such~that

\begin{iii}
\item[$\bullet$ ]
the~field extension\/
$\bbQ(\sqrt{D})/\bbQ$
splitting the double-six is~unramified
at\/~$p$,
\item[$\bullet$ ]
the~reduction\/
$S_p$
is geometrically irreducible and
no\/~\mbox{$\bbQ_p$-rational}
point
on\/~$S$
reduces to a singularity
of\/~$S_p$,
\end{iii}

\noindent
the~value
of\/~$\ev_p (\underline{c}, x)$
is independent
of\/~$x \in S(\bbQ_p)$.\medskip

\noindent
{\bf Proof.}
{\em
If~$p$
splits in the quadratic
extension~$\bbQ(\sqrt{D})$
then the assertion is true,~trivially. Thus,~we may assume that
$p$~remains
prime
in~$\bbQ(\sqrt{D})$.
The~requirement that
$\smash{z \in \bbQ_p^*}$
is a norm from
$\bbQ_p(\sqrt{D})$
then simply means that
$\nu_p(z)$~is~even.

Further,~we may restrict our considerations to
points~$x \in S(\bbQ_p)$
which are not contained in any of the 27~lines
on~$S$.
Indeed,~the local evaluation map is
\mbox{$p$-adically}
continuous and the complement of the 27~lines is dense
in~$S(\bbQ_p)$
according to Hensel's~lemma. In~particular, we may work
with~$F_{30}/F_{15}^2$~itself.

By~assumption, we have a
model~$\,\calS$
of~$S$
over~$\bbZ_p$
such that the special fiber
of~$\smash{\,\calS \times_{\Spec \bbZ_p} \Spec \calO_{\bbQ_p(\sqrt{D})}}$
is~irreducible. We~delete its singularities to obtain a
model~$\underline\calS$,
smooth
over~$\smash{\calO_{\bbQ_p(\sqrt{D})}}$.
According~to the last assumption,
every~$x \in S(\bbQ_p)$
determines a unique
extension~$\smash{\underline{x} \in \underline\calS(\calO_{\bbQ_p(\sqrt{D})})}$.

It~will suffice to construct a Zariski neighbourhood
of~$\underline{x}$
such that
$\ev_p (\underline{c}, \,.\; )$
is~constant. We~have, on the geometric generic~fiber,
$$\div(F_{30}/F_{15}^2) = 5E + 5G - 4F \, .$$
Here,~the divisors
$E := E_1 + \ldots + E_6$,
$G := G_1 + \ldots + G_6$,
and~$F := F_{12} + \ldots + F_{56}$
are
$\Gal(\overline\bbQ/\bbQ(\sqrt{D}))$-invariant,
and, therefore, defined
over~$S \times_{\Spec \bbQ_p} \Spec \bbQ_p (\sqrt{D})$.
$\underline\calS$~is
a regular model of that~variety. Hence,~every divisor
on~$\underline\calS$
is locally~principal. This~yields, in a Zariski
neighbourhood~$\calU_{\underline{x}}$,
$$F_{30}/F_{15}^2 = C p^k e^5 g^5 / f^4$$
for
$e$,
$g$,
and~$f$
rational functions corresponding to the divisors
$E$,
$G$,
and~$F$,
respectively,
$k \in \bbZ$,
and a
certain~$C \in \Gamma(\calU_{\underline{x}}, \calO^*_{\calU_{\underline{x}}})$.
Note~that we get by with one power
of~$p$
since the special fiber is~irreducible.

The~scheme
$\underline\calS$
is acted upon by the
conjugation~$\sigma \in \Gal(\bbQ(\sqrt{D})/\bbQ)$.
Restricting~to an open subscheme, if necessary, we may assume
that~$\smash{\calU_{\underline{x}}}$
is invariant
under~$\sigma$.
The~two sixers are interchanged
by~$\sigma$.
Consequently,
$$\sigma(e) = c p^l g, \quad \sigma(f) = c^\prime p^{l^\prime} f$$
for
$l, l^\prime \in \bbZ$
and regular
functions~$c, c^\prime$,
invertible
on~$\calU_{\underline{x}}$.
This~yields
$$F_{30}/F_{15}^2 = C c^\prime {}^2 C^{-5} p^{k+2l^\prime-5l} N(e^5/f^2) \, .$$
For~$x_0 \in S(\bbQ_p)$
specializing
to~$\calU_{\underline{x}}$,
we therefore see
$$\nu_p \big( (F_{30}/F_{15}^2) (x_0) \big) \equiv k+l \pmod 2 \, .$$
In~particular, the local evaluation map
$\ev_p (\underline{c}, \,.\; ) \colon S(\bbQ_p) \to \bbQ/\bbZ$
is constant on the set of all points specializing
to~$\calU_{\underline{x}}$.
As~the
point~$x \in S(\bbQ_p)$
defining the open
subset~$\calU_{\underline{x}}$
is arbitrary and the special fiber
$\calS_p$
is irreducible, this implies the~assertion.
}
\eop
\end{prop}

\begin{rems}
\begin{abc}
\item
Assuming~resolution of singularities in unequal characteristic, there is a proper
model~$\calS$
of~$S$
being a regular~scheme. Then,~for
$p$~a
prime unramified
in~$\bbQ(\sqrt{D})$,
the evaluation
$\ev_p (\underline{c}, x)$
depends only on the component
of~$\calS \times_{\Spec \bbZ} \Spec \bbF_{\!p^2}$,
the
point~$x$
specializes~to.
If~$p$
is ramified
and~$p \neq 2$
then we have at least that
$\ev_p (\underline{c}, x)$
is determined by the reduction
of~$x$
modulo~$p$.
\item
As one might expect from the proof given, Proposition~\ref{unverzw} is true in more~generality. The~reader might consult~\cite[Theorem~1]{B}.
\end{abc}
\end{rems}

\section{%Cubic surfaces obtained by
Explicit Galois~descent}

\begin{ttt}
Recall~that in~\cite{EJ}, we described a method to construct non-singular cubic surfaces
over~$\bbQ$
with a Galois invariant double-six. The~idea was to start with cubic surfaces in hexahedral~form. For~these, we developed an explicit version of Galois~descent.
\end{ttt}

\begin{ttt}
More~concretely, given a starting polynomial
$f \in \bbQ[T]$
of degree six without multiple zeroes, we construct a cubic
surface~$S_{(a_0,\ldots,a_5)}$
over~$\bbQ$
such~that
$$S_{(a_0,\ldots,a_5)} \times_{\Spec \bbQ} \Spec \overline\bbQ$$
is isomorphic to the
surface~$S^{(a_0,\ldots,a_5)}$
in~$\bP^5$
given~by
\begin{eqnarray*}
\phantom{a_0} X_0^3 + \phantom{a_1} X_1^3 + \phantom{a_2} X_2^3 + \phantom{a_3} X_3^3 + \phantom{a_4} X_4^3 + \phantom{a_5} X_5^3 & = & 0 \,
, \\
\phantom{a_0}\hsmash{X_0}\phantom{X_0^3} + \phantom{a_1}\hsmash{X_1}
\phantom{X_1^3} + \phantom{a_2}\hsmash{X_2}\phantom{X_2^3} + \phantom{a_3}\hsmash{X_3}\phantom{X_3^3} + \phantom{a_4}\hsmash{X_4}\phantom{X_4^3} + \phantom{a_5}\hsmash{X_5}\phantom{X_5^3} & = & 0 \, , \\
a_0 \hsmash{X_0}\phantom{X_0^3} + a_1 \hsmash{X_1}\phantom{X_1^3} + a_2 \hsmash{X_2}\phantom{X_2^3} + a_3 \hsmash{X_3}\phantom{X_3^3} + a_4 \hsmash{X_4}\phantom{X_4^3} + a_5 \hsmash{X_5}\phantom{X_5^3} & = & 0 \, .
\end{eqnarray*}
Here,~$a_0,\ldots,a_5 \in \overline\bbQ$
are the zeroes
of~$f$.

The~operation of an~element
$\sigma \in \Gal(\overline\bbQ/\bbQ)$
on~$S_{(a_0,\ldots,a_5)} \times_{\Spec \bbQ} \Spec \overline\bbQ$
goes over into the automorphism
$\smash{\pi_\sigma \!\circ\! t_\sigma \colon S^{(a_0, \ldots, a_5)} \to S^{(a_0, \ldots, a_5)} \, .}$
Here,~$\pi_\sigma$
permutes the coordinates according to the
rule~$a_{\pi_\sigma (i)} = \sigma(a_i)$
while
$t_\sigma$~is
the naive operation
of~$\sigma$
on~$S_{(a_0,\ldots,a_5)}$
as a morphism of schemes twisted
by~$\sigma$.
\end{ttt}

\begin{rems}
\begin{iii}
\item
More~details on the theory are given in~\cite[Theorem~6.6]{EJ}.
\item
On~$S_{(a_0,\ldots,a_5)}$,
there are the 15~obvious~lines given~by
$$X_{i_0} + X_{i_1} = X_{i_2} + X_{i_3} = X_{i_4} + X_{i_5} = 0$$
for~$\{i_0, i_1, i_2, i_3, i_4, i_5\} = \{0, 1, 2, 3, 4, 5\}$.
They~clearly form a Galois invariant~set. The~complement is a double-six. Correspondingly,~there are the 15~obvious tritangent~planes given
by~$X_i + X_j = 0$
for~$i \neq j$.

There~are formulas for the 30~non-obvious tritangent~planes,~too~\cite[Proposition~7.1.ii)]{EJ}. What~is important is that they are defined
over~$\bbQ(a_0, \ldots, a_5, \sqrt{d_4})$~for
$$\textstyle d_4(a_0, \ldots, a_5) := \sigma_2^2 - 4 \sigma_4 + \sigma_1 (2\sigma_3 - \frac32 \sigma_1\sigma_2 + \frac5{16} \sigma_1^3)$$
the Coble~quartic.
Here,~$\sigma_i$
is the
$i$-th
elementary symmetric function
in~$a_0, \ldots, a_5$.

Further,~an element
$\sigma \in \Gal(\overline\bbQ/\bbQ)$
flips the double-six if and only if it defines the conjugation
of~$\bbQ(\sqrt{D})$
for~$D := d_4 \!\cdot\! \Delta$,
the second factor denoting the discriminant
of~$a_0, \ldots, a_5$~\cite[Proposition~7.4]{EJ}.
\item
The~smooth manifold
$S(\bbR)$
has two connected components if and only if exactly four of the
$a_0, \ldots, a_5$
are real
and~$d_4(a_0, \ldots, a_5) > 0$.
Otherwise,~$S(\bbR)$
is~connected \cite[Corollary~8.4]{EJ}.
\item
The~descent
variety~$S_{(a_0,\ldots,a_5)}$
may easily be computed completely~explicitly. In~fact, \cite[Algorithm~6.7]{EJ} yields a quaternary cubic form with 20~rational~coefficients.
\end{iii}
\end{rems}

\begin{ttt}
\label{strat}
Using~the criteria provided in section~\ref{fuenf}, we have the following strategy to compute the Brauer-Manin obstruction
on~$S_{(a_0,\ldots,a_5)}$.\medskip\smallskip\pagebreak[3]

\noindent
{\bf Strategy}
(to explicitly compute the Brauer-Manin~obstruction on~$S_{(a_0,\ldots,a_5)}$).

\begin{iii}
\item
Compute~$D := d_4(a_0,\ldots,a_5) \!\cdot \! \Delta(a_0,\ldots,a_5)$.
Determine~the
list~$L_1$
of all primes at which
$\bbQ(\sqrt{D})$
is~ramified.
\item
By~a Gr\"obner basis calculation, determine all the primes
outside~$L_1$
at which
$S_{(a_0,\ldots,a_5)}$
has bad~reduction. Write~them into a
list~$L_2$.
\item
From~$L_2$,
delete all primes which split
in~$\bbQ(\sqrt{D})$.
Further,~erase all those primes
from~$L_2$
for which the singular points on the reduction
modulo~$p$
are not defined
over~$\bbF_{\!p}$
or do not lift
to~$S_{(a_0,\ldots,a_5)} \times_{\Spec \bbZ} \Spec \bbZ/p^k\bbZ$
for
$k$~large.
\item
Put~$L := L_1 \cup L_2$.
If~$D < 0$
and
$S_{(a_0,\ldots,a_5)}(\bbR)$
is not connected then the infinite place has to be added to this list of critical~primes.
\item
Delete all the primes
from~$L$
for which the
local~$H^1$-criterion
works~successfully.
\item
Scale~the
form~$F_{30}$
by a constant factor such that the local evaluation maps are zero for all the primes
outside~$L$.
\item
For~the
primes~$p$
that remained
in~$L$,
the
form~$F_{30}$
has to be~evaluated. For~that, cover
$S(\bbQ_p)$
by finitely many open subsets which are sufficiently small to ensure that the first
\mbox{$p$-adic}
digit
of~$F_{30}$
does not~change.
If~$p=2$
then the first three digits have to be taken into~account.

In~the case that we have a Galois invariant triple of azygetic double-sixes, the last step has to be executed three times, once for each of the corresponding forms of
type~$F_{30}$.
\end{iii}
\end{ttt}

\section{Application: Manin's conjecture% for the surfaces constructed
}

\begin{ttt}
Recall~that a conjecture, due to Yu.~I.~Manin, asserts that the number of 
\mbox{$\bbQ$-ra}\-tio\-nal
points of anticanonical height
$\leq\! B$
on a Fano
variety~$S$
is asymptotically equal to
$\tau B \log^{\rk \Pic (S) - 1} B$,
for~$B\to\infty$.
Further,~the coefficient
$\tau \in \bbR$
is conjectured to be the Tamagawa-type
number~$\tau (S)$
introduced by E.~Peyre in~\cite{Pe}. In~the particular case of a cubic surface, the anticanonical height is the same as the naive~height.
\end{ttt}

\begin{ttt}
E.~Peyre's Tamagawa-type~number is defined in \cite[Definition~2.4]{PT}~as
$$\smash{\tau(S) := \alpha(S) \!\cdot\! \beta(S) \cdot \lim_{s\to1} \, (s-1)^t L(s,\chi_{\Pic(S_{\overline\bbQ})}) \cdot \tau_H \!\big( S(\bbA_\bbQ)^{\rm Br} \big)}$$
for~$t = \rk \Pic(S)$.

Here,~the~factor
$\beta (S)$
is simply defined as
$\beta (S) := \# H^1 \big(\! \Gal(\overline\bbQ/\bbQ), \Pic(S_{\overline\bbQ}) \big)$.
$\alpha (S)$~is
given as follows~\cite[D\'efinition~2.4]{Pe}.
Let~$\Lambda_{\rm eff} (S) \subset \Pic(S) \otimes_\bbZ \bbR$
be the cone generated by the effective~divisors.
Identify~$\Pic(S) \otimes_\bbZ \bbR$
with~$\bbR^t$
via a mapping induced by an
isomorphism~$\smash{\Pic(S) \stackrel{\cong}{\longrightarrow} \bbZ^t}$.
Consider~the dual
cone~$\Lambda_{\rm eff}^\vee (S) \subset (\bbR^t)^\vee$.
Then,~$\alpha (S) := t \cdot \vol \, \{\, x \in \Lambda_{\rm eff}^\vee \mid \langle x, -K \rangle \leq 1 \,\}$.

$L(\;\cdot\;, \chi_{\smash{\Pic(S_{\overline\bbQ})}})$~denotes
the Artin
$L$-function
of the
$\Gal(\overline\bbQ/\bbQ)$-rep\-re\-sen\-tation
$\Pic(\smash{S_{\overline\bbQ}}) \otimes_\bbZ \bbC$
which contains the trivial representation
$t$~times
as a direct~summand. Therefore,
$L(s,\chi_{\Pic(S_{\overline\bbQ})}) = \zeta(s)^t \cdot L(s,\chi_P)$
and
$$\lim_{s\to1} \, (s-1)^t L(s,\chi_{\Pic(S_{\overline\bbQ})}) = L(1,\chi_P)$$
where
$\zeta$~denotes
the Riemann zeta function and
$P$~is
a representation which does not contain trivial~components. \cite[Corollary~11.5 and~Corollary~11.4]{Mu} show that
$L(s,\chi_P)$
has neither a pole nor a zero
at~$s=1$.

Finally,~$\tau_H$
is the {\em Tamagawa measure\/} on the set
$S(\bbA_\bbQ)$
of adelic points
on~$S$
and
$S(\bbA_\bbQ)^{\rm Br} \subseteq S(\bbA_\bbQ)$~denotes
the part which is not affected by the Brauer-Manin~obstruction.
\end{ttt}

\begin{ttt}
As~$S$
is projective, we~have
$\smash{S(\bbA_\bbQ) = \prod_{\nu \in \Val (\bbQ)} S(\bbQ_\nu)}$.
Then,~the Tamagawa measure
$\tau_H$
is defined to be the product measure
$\smash{\tau_H := \prod_{\nu \in \Val (\bbQ)} \tau_\nu}$.

Here,~for a prime
number~$p$,
the local
measure~$\tau_p$
on~$S(\bbQ_p)$
is given as~follows.
Let~$a \in S(\bbZ/p^k\bbZ)$
and
put~$\smash{\mathfrak{U}_{\kern0.4mm a}^{(k)}} := \{\, x \in S(\bbQ_p) \mid x \equiv a \pmodulo{p^k} \,\}$.
Then,
$$\tau_p (\mathfrak{U}_a^{(k)}) := \det (1 - p^{-1} \Frob_p \mid \Pic (S_{\overline\bbQ})^{I_p}) \cdot\! \lim_{m\to\infty} \!\!\!\frac{\#\{\, y \in S (\bbZ/p^m\bbZ) \mid y \equiv a \pmodulo{p^k} \,\}}{p^{m \dim S}} \, .$$
$\Pic (S_{\overline\bbQ})^{I_p}$~denotes
the fixed module under the inertia~group.

$\tau_\infty$
is described in~\cite[Lemme~5.4.7]{Pe}. In~the case of a cubic surface,
defined by the equation
$f = 0$,
this yields
$$\tau_\infty (U) = \frac12 \!\!\!\!\!\!\int\limits_{\substack{CU \\ |x_0|, \,\ldots\, , |x_3| \leq 1}} \!\!\!\!\!\!\!\! \omega_{\rm Leray}$$
for~$U \subset S (\bbR)$.
Here,~$\omega_{\rm Leray}$
is the {\em Leray measure\/} on the
cone~$CS (\bbR)$.
It~is related to the usual hypersurface measure by the formula\/
$\omega_{\rm Leray} = \smash{\frac1{\|\grad f\|}} \, \omega_{\rm hyp}$.
\end{ttt}

\begin{ttt}
Using~\cite[Algorithm~6.7]{EJ}, we constructed many examples of smooth cubic surfaces
over~$\bbQ$
with a Galois~invariant double-six. For~each of them, one may apply Strategy~\ref{strat} to compute the effect of the Brauer-Manin~obstruction. Then,~the method described in~\cite{EJ1} applies for the computation of Peyre's~constant.

From~the ample supply, the examples below were chosen in the hope that they indicate the main~phenomena. The~Brauer-Manin~obstruction may work at many primes simultaneously but examples where few primes are involved are the most~interesting. We~show that the fraction of the Tamagawa measure excluded by the obstruction can vary~greatly. We~also show that there may be an obstruction at the infinite~prime.
\end{ttt}

%\begin{ex}
%The~starting polynomial
%$$f = T^6 - 30T^4 + 20T^3 - 90T^2 - 1344T + 3970$$
%yields the cubic
%surface~$S$
%given by the~equation
%\begin{align*}
%& {-4}x^3 + 14x^2y - x^2z - 7x^2w + 14xy^2 + 2xyz + 8xyw - 4xzw \\
%& \hspace{2cm} - 8xw^2 - 9y^3 - 10y^2z + 8y^2w + yz^2 - 4yzw - 3yw^2 - 3zw^2 - 9w^3 = 0 \, .
%\end{align*}
%$S$~has
%bad reduction 
%at~$2$,
%$3$,
%$5$,
%$7$,
%$11$, 
%and~$31$.
%The~Galois group operating on the 27~lines
%is~$S_6 \times \bbZ/2\bbZ$
%and the orbit structure
%is~$[12, 15]$.
%The~quadratic field splitting the double-six
%is~$\bbQ(\sqrt{14})$.
%
%The~primes
%$5$,
%$11$,
%and~$31$
%split
%in~$\bbQ(\sqrt{14})$.
%Further,~the local evaluation map
%at~$7$
%turns out to be~constant. Hence,~the Brauer-Manin~obstruction works only at the primes
%$2$
%and~$3$.
%At~$2$,
%the local evaluation map decomposes
%$S(\bbQ_2)$
%into two sets of measures
%$\frac12$
%and~$\frac34$,
%respectively.
%At~$3$,
%the corresponding measures are
%$\frac{52}{243}$
%and~$\frac{214}{243}$.
%An~easy calculation
%shows~$\tau_H (S(\bbA_\bbQ)^{\rm Br}) = \frac{373}{665} \tau_H (S(\bbA_\bbQ))$.
%
%Using~this, for Peyre's constant, we find
%$\tau(S) \approx 2.7765$.
%Up~to a search bound
%of~$4000$,
%there are actually
%$10\,608$
%$\bbQ$-rational
%points in comparison with a prediction
%of~$11\,106$.
%\end{ex}

\begin{ex}
The~polynomial
$$f := T^6 - 390T^4 - 10\,180T^3 + 10\,800T^2 + 2\,164\,296T + 13\,361\,180 \in \bbQ[T]$$
yields the cubic
surface~$S$
given by the~equation
\begin{align*}
& {-x^2z} - x^2w - 3xy^2 + xz^2 + 14xzw + 8xw^2 - 2y^3 - 11y^2z \\
& \hspace{3cm} + y^2w + 4yz^2 + 4yzw + 10yw^2 + 4z^3 - 11z^2w + 9zw^2 - 6w^3 = 0 \, .
\end{align*}
$S$~has
bad reduction 
at~$2$,
$3$,
$5$,
$11$, 
and~$9\,265\,613\,761$.
The~Galois group operating on the 27~lines
is~$S_6$
acting in such a way that the orbit structure
is~$[12, 15]$.
Therefore,~we have
$\smash{H^1(\Gal(\overline{\bbQ}/\bbQ), \Pic(S_{\overline\bbQ})) \cong \bbZ/2\bbZ}$.
The~quadratic field splitting the double-six
is~$\bbQ(\sqrt{10})$.

The~primes
$3$
and~$9\,265\,613\,761$
split
in~$\bbQ(\sqrt{10})$.
The~local
$H^1$-criterion
excludes the
prime~$5$.
Further,~it turns out that the local evaluation map
at~$11$
is~constant. Hence,~the Brauer-Manin~obstruction works only at the
prime~$2$.
From~the whole
of~$S(\bbQ_2)$
which is of
measure~$4$
only a subset of
measure~$\frac94$
is~allowed.

Using~this, for Peyre's constant, we find
$\tau(S) \approx 1.7005$.
There~are actually
$6641$
$\bbQ$-rational
points of height at
most~$4000$
in comparison with a prediction
of~$6802$.
\end{ex}

\begin{ex}
The~polynomial
$$f := T^6 + 60T^4 - 40T^3 - 900T^2 + 15\,072T - 27\,860 \in \bbQ[T]$$
yields the cubic
surface~$S$
given by the~equation
\begin{align*}
& 5x^3 - 9x^2y + x^2z + 6x^2w + 3xy^2 + xyz + 6xyw - 2xz^2 \\
& \hspace{3cm} - 4xzw - y^3 - 3y^2z +  2yz^2 + 2yzw + 4z^3 + 2z^2w + 2zw^2 = 0 \, .
\end{align*}
$S$~has
bad reduction 
at~$2$,
$3$,
$5$, 
and~$73$.
The~Galois group operating on the 27~lines
is~$A_6 \times \bbZ/2\bbZ$
and the orbit structure
is~$[12, 15]$.
The~quadratic
field~$\bbQ(\sqrt{2})$
splits the double-six.

The~prime~$73$
splits
in~$\bbQ(\sqrt{2})$.
Further,~the local
$H^1$-criterion
excludes the
prime~$5$.
Hence,~the Brauer-Manin~obstruction works only at the primes
$2$
and~$3$.
At~$2$,
the local evaluation map decomposes
$S(\bbQ_2)$
into two sets of measures
$1$
and~$\frac12$,
respectively.
At~$3$,
the corresponding measures are
$\frac{70}{81}$
and~$\frac{28}{81}$.
An~easy calculation shows
$\tau_H (S(\bbA_\bbQ)^{\rm Br}) = \frac47 \tau_H (S(\bbA_\bbQ))$.

Using~this, for Peyre's constant, we find
$\tau(S) \approx 5.0879$.
Up~to a search bound
of~$4000$,
there are actually
$19\,302$
$\bbQ$-rational
points in comparison with a prediction
of~$20\,352$.
\end{ex}

\begin{ex}
The~polynomial
$f := T(T^5-5T-2) \in \bbQ[T]$
yields the cubic
surface~$S$
given by the~equation
\begin{align*}
& 2x^3 + x^2y - 4x^2z - x^2w + 2xy^2 + 2xyz + 2xyw - 2xz^2 - 4xzw \\
& \hspace{1.8cm} - 2xw^2 + 2y^2z - y^2w + yz^2 + 2yzw - 5yw^2 - 3z^2w + 6zw^2 + 9w^3 = 0 \, .
\end{align*}
$S$~has
bad reduction 
at~$2$,
$3$,
and~$5$.
Further,~$S(\bbR)$
consists of two connected~components. The~Galois group operating on the 27~lines is isomorphic
to~$S_5 \times \bbZ/2\bbZ$
and the orbit structure
is~$[12, 15]$.
$\bbQ(\sqrt{-15})$~is
the field splitting the double-six.

The~prime~$2$
splits
in~$\bbQ(\sqrt{-15})$.
Further,~the local
$H^1$-criterion
excludes the
prime~$5$.
Hence,~the Brauer-Manin~obstruction works only
at~$3$
and the infinite~prime.
At~$3$,
the local evaluation map decomposes
$S(\bbQ_2)$
into two sets of measures
$\frac23$
and~$\frac49$,
respectively.
At~the infinite prime, the corresponding measures are approximately
$1.9179$
and~$1.1673$.
An~easy calculation
shows~$\tau_H (S(\bbA_\bbQ)^{\rm Br}) \approx 0.5243 \!\cdot\! \tau_H (S(\bbA_\bbQ))$.

Using~this, for Peyre's constant, we find
$\tau(S) \approx 3.7217$.
Up~to a search bound
of~$4000$,
there are actually
$14\,249$
$\bbQ$-rational
points in comparison with a prediction
of~$14\,887$.
\end{ex}

\begin{ex}
The~polynomial
$$f := T(T^5 - 60T^3 - 90T^2 + 675T + 810) \in \bbQ[T]$$
yields the cubic
surface~$S$
given by the~equation
$$3x^3 + 2x^2z + xy^2 - 2xyz - 2xyw - xzw + 2xw^2 - yzw - yw^2 - z^3 + z^2w = 0 \, .$$
$S$~has
bad reduction 
at~$2$,
$3$,
and~$5$.
The~Galois group operating on the 27~lines is isomorphic
to~$S_5$
and the orbit structure
is~$[12, 15]$.
The~quadratic
field~$\bbQ(\sqrt{-3})$
splits the double-six.

The local
$H^1$-criterion
excludes the
prime~$5$.
Further,~the local evaluation maps turn out to be constant
on~$S(\bbQ_2)$
and~$S(\bbR)$.
At~the real prime, the reason is simply that
$S(\bbR)$~is~connected.
Consequently,~the Brauer-Manin~obstruction works only at the
prime~$3$.
From~the whole
of~$S(\bbQ_3)$
measuring~$\frac23$
a subset of
measure~$\frac49$
is~allowed.

Using~this, for Peyre's constant, we find
$\tau(S) \approx 2.2647$.
Up~to a search bound
of~$4000$,
there are actually
$8886$
$\bbQ$-rational
points in comparison with a prediction
of~$9059$.
\end{ex}

\begin{ex}
The~polynomial
$f := T(T^5+20T+16) \in \bbQ[T]$
yields the cubic
surface~$S$
given by the~equation
\begin{align*}
& {-3}x^3 - 7x^2y - 4x^2z + 5x^2w + 4xy^2 + 10xyz - 4xyw - 2xz^2 \\
& \hspace{3cm} + 2xzw + xw^2 - 4y^2z + yz^2 - 4yzw - 16yw^2 + z^2w - 5zw^2 = 0 \, .
\end{align*}
$S$~has
bad reduction 
at~$2$
and~$5$.
The~Galois group operating on the 27~lines is isomorphic
to~$A_5 \times \bbZ/2\bbZ$.
The~orbit structure
is~$[12, 15]$.
The~quadratic
field~$\bbQ(\sqrt{-5})$
splits the double-six.

The local
$H^1$-criterion
excludes the
prime~$5$.
At~the infinite prime, the local evaluation map is constant since
$S(\bbR)$
is~connected. Hence,~the Brauer-Manin~obstruction works only at the
prime~$2$.
It~allows a subset of
measure~$\frac{17}{16}$
out of
$S(\bbQ_2)$
measuring~$\frac54$.

Using~this, for Peyre's constant, we find
$\tau(S) \approx 2.4545$.
Up~to a search bound
of~$4000$,
there are actually
$9736$
$\bbQ$-rational
points in comparison with a prediction
of~$9818$.
\end{ex}

\begin{ex}
The~polynomial
$$f := T^6 - 456T^4 - 904T^3 + 102\,609T^2 + 1\,041\,060T + 2\,935\,300 \in \bbQ[T]$$
yields the cubic
surface~$S$
given by the~equation
\begin{align*}
& {-2}x^3 + 3x^2z + 9x^2w - 4xy^2 - 8xyz - 10xzw + 4xw^2 - 4y^3 - 3y^2z \\
& \hspace{3cm} - 4y^2w - 2yz^2 - 2yzw + 8yw^2 - z^3 + z^2w - 6zw^2 - 2w^3 = 0 \, .
\end{align*}
$S$~has
bad reduction 
at~$2$,
$3$,
$5$,
$31$,
and~$11\,071$.
The~Galois group operating on the 27~lines is isomorphic
to~$(S_3 \times S_3) \ltimes \bbZ/2\bbZ$
of order~72. The~orbit structure
is~$[6, 6, 6, 9]$.
There~is a triple of Galois~invariant double-sixes. Therefore,~we have
that~$\smash{H^1(\Gal(\overline{\bbQ}/\bbQ), \Pic(S_{\overline\bbQ})) \cong \bbZ/2\bbZ \times \bbZ/2\bbZ}$.
The~quadratic field splitting the double-sixes
is~$\bbQ(\sqrt{2})$.

The~primes
$31$
and~$11\,071$
split
in~$\bbQ(\sqrt{2})$.
Further,~the local evaluation maps turn out to be constant
on~$S(\bbQ_3)$.
Consequently,~the Brauer-Manin~obstruction works only at the primes
$2$
and~$5$.

The~local evaluation maps decompose
$S(\bbQ_2)$
into four sets of measures
$\frac7{16}$,
$\frac5{16}$,
$\frac14$,
and~$\frac14$,~respectively.
At~the
prime~$5$,
the corresponding measures are
$\frac{516}{625}$,
$0$,
$\frac{96}{625}$,
and~$0$.
Observe,~for one of the three non-zero Brauer classes, the local evaluation map is constant
on~$S(\bbQ_5)$.

A~simple calculation shows
that~$\tau_H (S(\bbA_\bbQ)^{\rm Br}) = \frac{111}{340} \tau_H (S(\bbA_\bbQ))$.
Using~this, for Peyre's constant, we find
$\tau(S) \approx 1.8532$.
Up~to a search bound
of~$4000$,
there are actually
$6994$
$\bbQ$-rational
points in comparison with a prediction
of~$7413$.
\end{ex}

\end{document}